\documentstyle[12pt]{article}
\makeatletter
\def\@sect#1#2#3#4#5#6[#7]#8{\ifnum #2>\c@secnumdepth
    \def\@svsec{}\else
    \refstepcounter{#1}\edef\@svsec{\csname the#1\endcsname.\hskip 1em }\fi
    \@tempskipa #5\relax
    \ifdim \@tempskipa>\z@
    \begingroup #6\relax
    \@hangfrom{\hskip #3\relax\@svsec}{\interlinepenalty \@M #8\par}
    \endgroup
    \csname #1mark\endcsname{#7}\addcontentsline
    {toc}{#1}{\ifnum #2>\c@secnumdepth \else
     \protect\numberline{\csname the#1\endcsname}\fi
           #7}\else
    \def\@svsechd{#6\hskip #3\@svsec #8\csname #1mark\endcsname
          {#7}\addcontentsline
          {toc}{#1}{\ifnum #2>\c@secnumdepth \else
     \protect\numberline{\csname the#1\endcsname}\fi
           #7}}\fi
     \@xsect{#5}}
\def\label#1{\@bsphack\if@filesw {\let\thepage\relax
   \xdef\@gtempa{\write\@auxout{\string
   \newlabel{#1}{{\thesection.\@currentlabel}{\thepage}}}}}\@gtempa
   \if@nobreak \ifvmode\nobreak\fi\fi\fi\@esphack}
\def\@eqnnum{(\thesection.\theequation)}
\def\section{\setcounter{equation}{0} \@startsection {section}{1}{\z@}{-3.5ex
   plus -1ex minus -.2ex}{2.3ex plus .2ex}{\Large\bf}}
\newcount\@minsofar
\newcount\@min
\newcount\@cite@temp
\def\@citex[#1]#2{%
\if@filesw \immediate \write \@auxout {\string \citation {#2}}\fi
\@tempcntb\m@ne \let\@h@ld\relax \def\@citea{}%
\@min\m@ne%
\@cite{%
  \@for \@citeb:=#2\do {\@ifundefined {b@\@citeb}%
    {\@h@ld\@citea\@tempcntb\m@ne{\bf ?}%
    \@warning {Citation `\@citeb ' on page \thepage \space undefined}}%
{\@minsofar\z@ \@for \@scan@cites:=#2\do {%
  \@ifundefined{b@\@scan@cites}%
    {\@cite@temp\m@ne}
    {\@cite@temp\number\csname b@\@scan@cites \endcsname \relax}%
\ifnum\@cite@temp > \@min
    \ifnum\@minsofar = \z@
      \@minsofar\number\@cite@temp
      \edef\@scan@copy{\@scan@cites}\else
    \ifnum\@cite@temp < \@minsofar
      \@minsofar\number\@cite@temp
      \edef\@scan@copy{\@scan@cites}\fi\fi\fi}\@tempcnta\@min
  \ifnum\@minsofar > \z@ 
    \advance\@tempcnta\@ne
    \@min\@minsofar
    \ifnum\@tempcnta=\@minsofar 
      \ifx\@h@ld\relax
        \edef \@h@ld{\@citea\csname b@\@scan@copy\endcsname}%
    \else \edef\@h@ld{\ifmmode{-}\else--\fi\csname b@\@scan@copy\endcsname}%
      \fi
    \else \@h@ld\@citea\csname b@\@scan@copy\endcsname
          \let\@h@ld\relax
  \fi 
\fi}%
\def\@citea{,\penalty\@highpenalty\,}}\@h@ld}{#1}}
\makeatother
\begin{document}
\def\t{\theta}
\def\T{\Theta}
\def\D{\Delta}
\def\w{\omega}
\def\ov{\overline}
\def\a{\alpha}
\def\b{\beta}
\def\g{\gamma}
\def\s{\sigma}
\def\l{\lambda}
\def\wt{\widetilde}
\def\di{\displaystyle}
\def\sn{\mbox{sn}}
\def\cd{\mbox{cd}}
\def\cn{\mbox{cn}}
\def\dn{\mbox{dn}}
\def\L{\Lambda_{n,{\bf{Q}}(q,t)}}
\def\be{\begin{equation}}
\def\ee{\end{equation}}
\def\beq{\begin{eqnarray}}
\def\eeq{\end{eqnarray}}
\def\oo{\circ}
\def\ds{\displaystyle}
\newtheorem{lemma}{Lemma}[section]
\newtheorem{theorem}[lemma]{Theorem}

\addtolength{\unitlength}{-0.5\unitlength}

\newsavebox{\square}\savebox{\square}(30,27){\begin{picture}(30,27)
\thicklines\multiput(10,0)(20,0){2}{\line(0,1){20}}
\multiput(10,0)(0,20){2}{\line(1,0){20}}\end{picture}}
\hspace{10cm}November, 1995.
\vspace{3cm}

\centerline{\bf\Large $A_2$ Macdonald polynomials:}
\vspace{0.5cm}
\centerline{\bf\Large  a separation of variables}
\vspace{1cm}
\centerline{
V.V. Mangazeev\footnote{
E-mail: vvm105@phys.anu.edu.au}}
\vspace{0.2 cm}
\centerline{\small Department of Theoretical Physics, RSPhysSE,}
\centerline{\small
Australian National University,}
\centerline{\small
 Canberra, ACT 0200, Australia}
\vspace{1cm}

\vspace{1cm}
\centerline{\bf Abstract}
In this paper we construct a discrete linear operator $K$ which
transforms $A_2$ Macdonald polynomials into the product of two
basic $\phantom{|}_3\phi_2$ hypergeometric series  with known
arguments. The action of the operator $K$ on power sums in two
variables can be reduced to a generalization of one particular
case of the Bailey's summation formula for a very-well-poised
$\phantom{|}_6\psi_6$ series. We also propose the conjecture
for a transformation of $\phantom{|}_6\psi_6$ series with
different arguments.

\newpage

\section{Introduction}

In the recent paper \cite{KS} V.B. Kuznetsov and E.K. Sklyanin
proposed a new integral representation for $A_2$ Jack polynomials
\cite{Ja} in terms
of $\phantom{|}_3F_2$ hypergeometric functions and constructed the integral
operator $M$ and its inverse $M^{-1}$ which separates variables in
$A_2$ Jack polynomials. Also they formulated several conjectures
about a structure of separated functions for general $A_n$ case.
In fact, their method of a separation
of variables originates (see recent review
\cite{SK}) from the Inverse Scattering Method \cite{FT}  and
$L$-operator formalism applied to Calogero-Sutherland system \cite{Ca,Su}.
It is known that in this case the corresponding classical $r$-matrix
depends on dynamical variables \cite{AT,SK1,BS} and does not satisfy
the usual classical Yang-Baxter equation. As a result we don't know
a consistent quantization procedure for such kind of integrable systems.
Nevertheless, using the classical $L$-operator authors of ref. \cite{KS}
guessed a quantum version of the equations on the integral
operator $M$ which permits a separation of variables for $A_2$
Jack polynomials.

It is well known that the quantum Calogero-Sutherland system has
a discrete $q$-finite analog \cite{MD}. Corresponding eigenfunctions are
called Macdonald $A_n$ polynomials which generalize Jack $A_n$
polynomials for finite $q$. In this paper we will show that
a separation variable procedure should work for Macdonald polynomials
as well. We will formulate some results concerning general $A_n$
Macdonald polynomials and construct explicitly for $A_2$ case
a discrete operator $K$ which apparently separates variables
in $A_2$ Macdonald polynomials (we believe that this is true, but
a rigorous proof  can be reduced to summations of some hypergeometric
$q$-series). The paper is organized as follows.

In Section 2 we remind some well known facts about Calogero-Sutherland
system. Section 3 contains a brief review of Macdonald symmetric
functions related to the root system $A_n$. In Section 4 we give
definitions related to basic hypergeometric series and two theorems
which show a connection between the spectrum of $A_{n-1}$ Macdonald
operators and $\phantom{|}_n\phi_{n-1}$ series.
In Section 5 we construct explicitly for $A_2$ case a discrete
operator $K$ which separates variables in  Macdonald polynomials and
formulate the conjecture for the action of this operator on power
sums in two variables. Section 6 contains a discussion of some
unsolved problems. In Appendix A we write out a set of difference
equations on the kernel of the operator $K$
which appear from separated equations of Section 4
for $A_2$ case. At last, Appendix B contains the calculation of
the normalization factor for the kernel of the operator $K$ and
several examples which support the conjecture from Section 5.

\section{Calogero-Sutherland model}

In this section we remind some known results related to
Calogero-Sutherland hamiltonian.

Consider $n$ particles on a circle interacting with a long range
potential \cite{Ca,Su}. We denote coordinates of the particles by
$x_i$, $i=1,\ldots,n$, $0\le x_i\le L$. Then a total momentum
of the system and hamiltonian which describe a dynamics of the particles
are given by:
\be
P=\sum_{i=1}^n{\ds1\over\ds i}{\ds d\over\ds dx_i},\label{K1}
\ee
\be
H=-\sum_{i=1}^n{\ds 1\over\ds2}{\ds d^2\over\ds dx_i^2}+g(g-1)
{\ds\pi^2\over\ds L^2}
\sum_{i<j}{\ds 1\over \ds \sin^2({ \pi\over L}(x_i-x_j))}. \label{K2}
\ee
As usual set $\theta_j=2\pi{\ds x_j\over\ds L}$ and $z_j=\exp(i\theta_j)$.
Then eigenfunctions of the equation $H\Psi=E\Psi$ have the following
structure
\be
\Psi(\theta)=\Delta^g(\theta)J(\theta),\label{K3}
\ee
\be
\Delta(\theta)=\prod_{i<j}\sin\Bigl({\ds\theta_i-\theta_j\over\ds2}\Bigl)
\label{K4}
\ee
and $J$ is a symmetric Laurent polynomial in $z_j$.

It is naturally to introduce the effective Hamiltonian
\be
\tilde H=\Delta^{-g}H\Delta^g, \label{K5}
\ee
which has polynomial eigenfunctions $J$ and can be rewritten
in the following form
\be
\tilde H=\sum_{j=1}^n\bigl(z_j{\ds\partial\over\ds\partial z_j}\bigl)^2+
g\sum_{i\ne j}{\ds z_i+z_j\over\ds z_i-z_j}
\bigl(z_i{\ds\partial\over\ds\partial z_i}-
z_j{\ds\partial\over\ds\partial z_j}\bigl).\label{K6}
\ee
Now we list some well known facts about operator $\tilde H$.
\begin{enumerate}
\item
$\tilde H$ can be included in a mutually commuting family of
differential operators $\tilde H_k$, $k=1,\ldots,n$ in variables $z_j$
\cite{Sek}.
\item
Simultaneous polynomial eigenfunctions $J$ of the
differential operators $\tilde H_k$ are known as Jack polynomials.
\item
This system has a discrete $q$-analog \cite{MD}. Corresponding commuting
discrete operators are called Macdonald operators. Their simultaneous
eigenfunctions are Macdonald polynomials which generalize
Jack polynomials for finite $q$.
\end{enumerate}

\section{Macdonald symmetric functions}

In this section we give a short list of definitions concerning
Macdonald polynomials \cite{MD} related to $A_{n-1}$ root system.
There are several ways to define Macdonald
symmetric functions. Keeping in mind a connection with
Calogero-Sutherland model we define Macdonald polynomials $P_\l(x;q,t)$
as simultaneous eigenfunctions of a commuting family of operators
$D_n^r$, $r=0,1,\ldots,n$ acting in the ring of symmetric functions
in $n$ variables.

As usual let $\L$ be the ring of symmetric functions in $n$ variables
$(x_1,\ldots,x_n)$ over field ${\bf{Q}}(q,t)$ of rational functions
in two independent indeterminates $q$ and $t$.

Define shift operators $T_{x_i}$ by

\begin{equation}
(T_{x_i}f)(x_1,\ldots,x_n)=f(x_1,\ldots,qx_i,\dots,x_n)   \label{M1}
\end{equation}
for any polynomial $f(x_1,\ldots,x_n)$.

Now let us introduce a family of mutually commuting
operators $D_n^r$, $r=0,1,\ldots,n$ as follows

\begin{equation}
D_n^r=\sum_I t^{r(r-1)/2}
\prod_{i\in I \atop j\not\in I}{tx_i-x_j\over x_i-x_j}\prod_{i\in I}
T_{x_i},        \label{M2}
\end{equation}
where the summation is over all $r$-element subsets $I$ of $(1,2,\ldots,n)$.

Define a partition $\l$ of the weight $|\l|$ as a
sequence of non-negative integers $(\l_1,\l_2,\ldots)$, $\l_1\ge\l_2\ge\ldots$
such that $|\l|=\sum_i \l_i < \infty$. The nonzero's $\l_i$
are called the parts of $\l$ and the number of parts is the length
$l(\l)$ of the partition $\l$.
A natural partial ordering for two partitions $\l, \mu$ can be defined
as follows
\begin{equation}
\l\ge\mu \Leftrightarrow |\l|=|\mu|\>\mbox{and}\>
\sum_{i=1}^r\l_i\ge\sum_{i=1}^r\mu_i\>
\mbox{for all r}\ge 1.    \label{M3}
\end{equation}

One of the bases in $\L$ is given by monomial symmetric functions $m_\l$,
$l(\l)\le n$:
\begin{equation}
m_\l=\sum_\a x^{\a(\l)},   \label{M4}
\end{equation}
where $x^\l=x_1^{\l_1}x_2^{\l_2}\ldots$ and
the summation is over all distinct permutations $\a$ of nonzero  parts of
the partition $\l$.

Further for each partition $\l$, $l(\l)\le n$ let coefficients
$d_n^i(\l)$ are defined by the following expansion
\begin{equation}
\prod_{i=1}^n(1+Xq^{\l_i}t^{n-i})=\sum_{i=0}^n X^i d_n^i(\l).   \label{M5}
\end{equation}

\begin{theorem}{\rm \cite{MD}}
For each partition $\l$, $l(\l)\le n$ there exists a unique symmetric
function $P_\l(x;q,t)\in\L$ which satisfies two following conditions:
\begin{eqnarray}
  1.
  &&
  P_{\lambda}=m_{\lambda} + \sum_{\mu<\lambda} u_{\lambda\mu}m_{\mu},
  \quad
  u_{\lambda\mu}\in {\bf Q}(q,t),\label{M6} \\
  2.
  &&
  D_n^r P_{\lambda}=d_n^r(\l)\> P_{\lambda}.\label{M7}
\end{eqnarray}
\end{theorem}

One can also define the scalar product $\langle.\>,\>.\rangle_{q,t}$ \cite{MD}
such that
\be
\langle P_\l,P_\mu\rangle_{q,t}=0,\quad \mbox{if}\quad\l\ne\mu,  \label{M8}
\ee
but we will not use it in this paper.

\section{Spectrum of Macdonald operators and basic hypergeometric series}

In this section we remind a definition of the basic hypergeometric series
\hspace{-0.2cm}$\phantom{a}_{r+1}\phi_{r}$ and show that coefficients
$d_n^i(\l)$ (see (\ref{M5})) naturally appear from separated equations for
the function \hspace{-0.2cm}$\phantom{a}_n\phi_{n-1}$.

First of all introduce the following notations
(see, for example, \cite{GR,An}):
\begin{equation}
(a;q)_\infty=\prod_{k=0}^\infty(1-aq^k), \label{S1}
\end{equation}
\be
(a;q)_n={(a;q)_\infty\over (aq^n;q)_\infty},\quad n\in{\bf Z},
\label{S2}
\ee
\be
(a_1,a_2,\ldots,a_r;q)_n=\prod_{i=1}^r(a_i;q)_n,
\label{S3}
\ee
where $|q|<1$, $a\in{\bf C}$, $r\ge0$.

As usual define the basic hypergeometric function
\hspace{-0.2cm}$\phantom{a}_{r+1}\phi_{r}$ as
\be
\phantom{a}_{r+1}\phi_{r}(a_1,\ldots,a_{r+1};b_1,\ldots,b_r;q,z)=
\sum_{n=0}^\infty{(a_1;q)_n(a_2;q)_n\ldots(a_{r+1};q)_n\over
(q;q)_n(b_1;q)_n\ldots(b_r;q)_n}z^n.   \label{S4}
\ee
Sometimes we will use a short
notation \hspace{-0.2cm}$\phantom{a}_{r+1}\phi_{r}([a];[b];q,z)$.

This function satisfies the difference equation
\be
\Bigl[(1-T_{z})\prod_{i=1}^r(1-q^{-1}b_iT_{z})-
z\prod_{q=1}^{r+1}(1-a_iT_{z})\Bigr]
\hspace{-0.1cm}\phantom{a}_{r+1}\phi_{r}([a];[b];q,z)=0.  \label{S5}
\ee

Now for each partition $\l$, $l(\l)\le n$ define two sequences
$[a]=(a_1,\ldots,a_{n})$, $[b]=(b_1,\ldots,b_{n-1})$
\be
a_i=q^{1+\l_n-\l_i}t^{i-1-n},\> i=1,\ldots,n ,\quad
b_j=ta_j,\>j=1,\ldots,n-1               \label{S6}
\ee
and introduce  functions $\varphi_\l(q,t;z)$, $\psi_\l(q,t;z)$:
\be
\varphi_\l(q,t;z)=z^{\l_n}
\hspace{-0.2cm}\phantom{a}_{n}\phi_{n-1}([a];[b];q,z)   \label{S7}
\ee
\be
\psi_\l(q,t;z)={(z,q)_\infty \over (qzt^{-n},q)_\infty}\>\varphi_\l(q,t;z).
\label{S8}
\ee

\noindent
Then
\begin{theorem}
The function $\psi_\l(q,t;z)$ is a polynomial in $z$.
\end{theorem}
Proof:

\noindent
Using (\ref{S6}) we note that
\be
{\ds (a_{i+1};q)_l\over\ds(b_i;q)_l}=
{\ds (b_i q^l;q)_{\l_i-\l_{i+1}}\over\ds(b_i;q)_{\l_i-\l_{i+1}}},\quad
i=1,\ldots,n-1.\label{S9}
\ee
Then
\beq
&
{\ds\psi_\l(q,t;z)=z^{\l_n}{(z,q)_\infty \over (qzt^{-n},q)_\infty}
\sum_{l=0}^\infty{\ds (q^{1+\l_n-\l_1}t^{-n};q)_l\over\ds (q;q)_l}
\prod_{i=1}^{n-1}\biggl[
{\ds (b_i q^l;q)_{\l_i-\l_{i+1}}\over\ds(b_i;q)_{\l_i-\l_{i+1}}}\biggl]z^l=}&
\nonumber\\
&
={\ds z^{\l_n}{(z,q)_\infty \over (qzt^{-n},q)_\infty}
\sum_{l=0}^\infty{\ds (q^{1+\l_n-\l_1}t^{-n};q)_l\over\ds (q;q)_l}
\sum_{k=0}^{\l_1-\l_n}\rho_k(t;q)(q^kz)^l}=&\nonumber\\
&
={\ds z^{\l_n}{(z,q)_\infty \over (qzt^{-n},q)_\infty}
\sum_{k=0}^{\l_1-\l_n}\rho_k(t;q)
{\ds (q^{1+\l_n-\l_1+k}t^{-n}z;q)_\infty\over\ds (q^kz;q)_\infty}
}=&\nonumber\\
&
={\ds z^{\l_n}
\sum_{k=0}^{\l_1-\l_n}\rho_k(t;q)
{\ds (z;q)_k\over\ds
(qzt^{-n};q)_{\l_n-\l_1+k}}
}={\ds\sum_{k=\l_n}^{\l_1}\tilde
\rho_k(t;q)z^k}.\quad\mbox{\it Q.E.D.}&\nonumber
\eeq

\noindent
In fact, this theorem naturally generalizes {\bf
Theorem 3} from ref. \cite{KS}.

And now we formulate the main result of this section:

\begin{theorem}
The function $\varphi_\l(q,t;z)$ satisfies the following difference
equation
\be
\Bigl[\sum_{i=0}^n(1-z(qt^{-1})^i)d_n^{n-i}(\l)(-T_{z})^i\bigr]
\,\varphi_\l(q,t;z)=0,     \label{S10}
\ee
where coefficients $d_n^i(\l)$ coincide with eigenvalues of Macdonald
operators $D_n^i$ for the eigenfunction $P_\l(x;q,t)$.
\end{theorem}
Proof:
Direct calculation. {\it Q.E.D.}

We will call equation (\ref{S10}) as the separated equation for $A_{n-1}$
Macdonald polynomials. This equation strongly supports (and generalizes)
the conjecture of ref. \cite{KS} concerning a structure of separated functions
for $A_{n-1}$ Jack polynomials.
In the next section we will construct explicitly a discrete linear
operator $K$ which provides a separation of variables for $A_2$
Macdonald polynomials.

\section{Separation of variables for $A_2$ Macdonald polynomials}.

Hereafter we will consider only the case $n=3$, which corresponds
to $A_2$ Macdonald polynomials.

Using (\ref{M2}) we can write out explicitly Macdonald operators for $n=3$:

\beq
&&D_3^0=1, \nonumber\\
&&D_3^1={(tx_1-x_2)(tx_1-x_3)\over(x_1-x_2)(x_1-x_3)}T_{x_1}+
(x_1\to x_2\to x_3),\nonumber\\
&&D_3^2=t{(tx_1-x_3)(tx_2-x_3)\over(x_1-x_3)(x_2-x_3)}T_{x_1}T_{x_2}+
(x_1\to x_2\to x_3),\nonumber\\
&&D_3^3=t^3T_{x_1}T_{x_2}T_{x_3}.\label{A1}
\eeq

First let us make the following changing of variables:
\be
x_1=z_1v,\quad x_2=z_2v,\quad x_3=v, \label{A2}
\ee
\be
T_{x_1}=T_{z_1},\quad T_{x_2}=T_{z_2},\quad T_{x_3}=T^{-1}_{z_1}
T^{-1}_{z_2}T_v.  \label{A3}
\ee
Then using a homogeneity of
$A_2$ Macdonald polynomials $P_\l(x_1,x_2,x_3;q,t)$
with respect to $x_i$ we have
\be
P_\l(x_1,x_2,x_3;q,t)=v^{\l_1+\l_2+\l_3}p_\l(z_1,z_2;q,t), \label{A4}
\ee
where $p_\l(z_1,z_2;q,t)$ is a nontrivial function of two new  variables
$z_1$ and $z_2$.

Now let us suppose an existence of a linear discrete operator $K$
\be
K:\> p_\l(z_1,z_2)\>\to\>(K\oo p_\l)(z_1,z_2)\label{A5}
\ee
such that after applying of the $K$ to each $p_\l(z_1,z_2)$ a result
will satisfy to separated equation (\ref{S10}) (at $n=3$) in each
variable $z_i$.
In fact, the action of this linear operator can be represented
in the form of the summation of
$p_\l$ with some discrete kernel $K$.

Then using a linearity of the operator $K$ let us rewrite
separated equations (\ref{S10}) as
\be
\bigl\{\sum_{l=0}^3(1-z_i(qt^{-1})^l)(-T_{z_i}
\begin{array}[t]{ccc}
)^l\bigr\}(K & \oo & d_3^{3-l}(\l)\\
\phantom{a}{\Uparrow} & & {\Downarrow}\\
D_3^{3-l\>*} & \Leftarrow & D_3^{3-l}\end{array}
p_\l)(z_1,z_2)=0.\label{A6}
\ee
Now we can replace eigenvalues $d_3^{3-l}(\l)$ of Macdonald polynomials
by corresponding Macdonald operators and using a conjugation with
respect to the convolution $K\oo p_\l$
to move the action
of Macdonald operators on the kernel of the operator $K$.
In fact, it will correspond to shifts of summation variables
on $+1$ or $-1$.
After that we can try to satisfy a set of local difference equations
on the kernel $K$. And if we solve these equations, probably
this linear operator $K$ will separate variables in $A_2$ Macdonald
polynomials.

But if summations are not infinite in both directions, some problems
with boundary terms can appear under shifts of summation variables.
However, it appears that all summations can be chosen from
$-\infty$ to $+\infty$ and for each fixed Macdonald polynomial
$p_\l$ we can choose parameter $t$ in such a way
that all series will be convergent.
So we will not care about boundary terms anymore.

To write out explicitly the action of the linear operator $K$ on polynomials
$p_\l(z_1,z_2)$ let us make the following substitution:
\be
z_1=\sigma\xi,\quad z_2=\sigma\xi^{-1}\label{A7}
\ee
\be
T_{z_1}=T_\sigma^{1/2}T_{\xi}^{1/2},\quad
T_{z_2}=T_\sigma^{1/2}T_{\xi}^{-1/2}.\label{A8}
\ee
It turns out that the action of the operator $K$
on the variable $\sigma$ is very simple:
\be
K\oo \sigma^l=t^{3l/2}\>\sigma^l,\label{A9}
\ee
but the action of $K$ on the variable $\xi$ is quite nontrivial.

So taking into account formulas (\ref{A7}-\ref{A9}) we can represent
the action of $K$ in the form of the following single sum:
\be
K:\> p_\l\>\to\>\sum_{l=-\infty}^{+\infty}
\> K(z_1,z_2;y_1,y_2)\>p_\l(y_1,y_2),\label{A10}
\ee
where
\be
z_1=t^{3/2}\sigma\xi,\quad z_2=t^{3/2}\sigma\xi^{-1}, \label{A11}
\ee
\be
y_1=\sigma\eta,\quad y_2=\sigma\eta^{-1},\quad
\eta=q^{s+l/2},\quad l\in Z \label{A12}
\ee
and $s\in[0,1/2]$. We will explain a necessity to introduce the
parameter $s$ later.

Further we will use one of the following equivalent forms for the
kernel $K$:
\be
K(z_1,z_2;y_1,y_2)=K(z;y)=K(\sigma,\xi;\eta)=
K(z_1,z_2;s,l), \label{A13}
\ee
where arguments satisfy (\ref{A11}-\ref{A12}).

Note that due to (\ref{A8}) we must summarize, in fact, over integer
and half-integer powers of $q$. But we can split the summation in (\ref{A10})
into two sums: with the kernel $K(z_1,z_2;s,l)$ and with
the kernel $K(z_1,z_2;s+1/2,l)$ and summarize over $l$ with $\eta=q^{s+l}$.

Substituting (\ref{A10}) into separated equation (\ref{A6}) and using
explicit form (\ref{A1}) of the Macdonald operators
one can obtain four linear difference equations
on the kernel $K$. These equations are given in Appendix A.
It turns out that they are a consequence of much more simple
relations:

\be
T_{z_i}T_{y_1}K(z;y)={\ds t^2(ty_1-1)(y_2-qy_1)(z_i-ty_2)\over
\ds(qy_1-t)(y_2-y_1)(qz_i-t^2y_2)}K(z;y),\label{UR1}
\ee
\be
T_{z_i}T_{y_2}K(z;y)={\ds t^2(ty_2-1)(y_1-qy_2)(z_i-ty_1)\over
\ds(qy_2-t)(y_1-y_2)(qz_i-t^2y_1)}K(z;y).\label{UR2}
\ee

\vspace{0.3cm}
Let us make several comments about equations (\ref{UR1}-\ref{UR2}).
If we consider the limit $q\to 1$, then it can be shown that these
equations will reduce to differential equations on the kernel of
the integral operator which separates variables for $A_2$ Jack
polynomials \cite{SK}. In fact, E.K. Sklyanin showed that one can
extract these equations (their classical version) from classical
$L$-operator related to Calogero-Sutherland system. We obtained
(\ref{UR1}-\ref{UR2}) explicitly from the separated equations
(\ref{A6}) for $A_2$ Macdonald polynomials. It would be interesting
to understand: how  to obtain  equations (\ref{UR1}-\ref{UR2})
from some quantum $q$-finite version of the classical $L$-operator
related to $3$-particle Calogero-Sutherland system.
For $n>3$ it is rather difficult to work with a generalization
of difference equations (\ref{Ap1}-\ref{Ap2}) (see Appendix A) because
of their very complicated structure.
But it is likely that equations (\ref{UR1}-\ref{UR2}) should have
a simple multiplicative generalization.

Now it is not difficult to solve equations (\ref{UR1}-\ref{UR2})
with respect to the kernel $K$.
The answer is
\be
K(\sigma,\xi;\eta)=
g(\sigma,\xi,s)
{\ds\eta(1-\eta^2)\over\ds(t\sqrt{t})^{s+l}}
{\ds \biggl(q{\eta\over\sigma t},q{\sigma\eta\over t},
q{\eta\xi\over\sqrt{t}},q{\eta\over\xi\sqrt{t}};q\biggl)_\infty\over
\ds\biggl(t\eta\sigma,{t\eta\over\sigma},\eta\xi\sqrt{t},
{\eta\sqrt{t}\over\xi};q\biggl)_\infty},\label{A16}
\ee
where
\be
{\ds g(q^{1/2}\sigma,q^{\pm 1/2}\xi,s)\over\ds g(\sigma,\xi,s)}=
{\ds q^{1/2}\over\ds t^{3/2}}.\label{A17}
\ee

Let us note that if we set the parameter $s$ in (\ref{A12}) as $s=0$
then the sum in RHS of (\ref{A10}) will be equal zero (it can be checked
explicitly using (\ref{A16})). So the function $g(\sigma,\xi,s)$
must have a pole at $s=0$. We will see that this is really the case.

It is enough to calculate the action of the operator $K$
on the power sums $z_1^n+z_2^n$.

Define the following function
\beq
&\Phi(s,\sigma,\xi,t;q,n)=
{\ds\sum_{l=-\infty}^{\infty}}
{\ds q^{s+l}(1-q^{2(s+l)})\over\ds t^{3(s+l)}}\times&\nonumber\\
&\times(q^{n(s+l)}+q^{-n(s+l)})
{\ds\Bigl({\ds q\sigma\over \ds t}q^{s+l},{\ds q\over\ds t\sigma}q^{s+l},
{\ds q\over\ds\xi\sqrt{t}}q^{s+l},{\ds q\xi\over\ds\sqrt{t}}q^{s+l}
\Bigl)_\infty
\over
{\ds\Bigl({\ds t\over \ds \sigma}q^{s+l},{t\sigma}q^{s+l},
{\xi\sqrt{t}}q^{s+l},{\ds \sqrt{t}\over\ds\xi}q^{s+l}
\Bigl)_\infty}},\label{A18}&
\eeq
where $s\in[0,1]$, $n=0,1,2,3,\ldots$

It is easy to see that all series in (\ref{A18}) are convergent
provided that $|q^{1-n}/t^3|<1$. So for any fixed $n$ we can always choose
parameter $t$ in such a way that sums in  (\ref{A18}) will be convergent
and after analytically continue a result for any complex $t$.

Then we have
\be
K\oo(z_1^n+z_2^n)=g(\sigma,\xi,s)t^{3n/2}\sigma^n(\Phi(s,\sigma,\xi,t;q,n)+
\Phi(s+{1\over2},\sigma,\xi,t;q,n))  \label{A19}
\ee
where $z_1$, $z_2$ satisfy (\ref{A7}).

It is naturally to choose the function $g(\sigma,\xi,s)$ as follows
\be
g(\sigma,\xi,s)={\ds2\over\ds\Phi(s,\sigma,\xi,t;q,0)+
\Phi(s+1/2,\sigma,\xi,t;q,0)}  \label{A20}
\ee
and as a consequence
\be
K\oo 1=1.\label{A21}
\ee

Using expression (\ref{B2}) (see Appendix B) for the function
$\Phi(s,\sigma,\xi,t;q,0)$ it is not difficult to show that the function
$g(\sigma,\xi,s)$ defined by (\ref{A20}) satisfies reccurence relations
(\ref{A17}).

Now let us formulate the following

\vspace{0.5cm}
{\bf Conjecture:}
{\it
\be
\Phi(s,\sigma,\xi,t;q,n)=\Phi(s,\sigma,\xi,t;q,0)P(\sigma,\xi,t;q,n),
\label{A22}
\ee
where $P(\sigma,\xi,t;q,n)$ is a Laurent polynomial in
$\sigma$, $\xi$.}

Three first nontrivial examples of polynomials $P(\sigma,\xi,t;q,n)$
are given in Appendix B.

{}From (\ref{A19}-\ref{A22}) it is easy to see that the action of the
operator $K$ on power sums does not depend on the parameter $s$.

So we immediately come to the following statement:

\vspace{0.3cm}
{\it If the conjecture (\ref{A22}) is valid for the arbitrary $n$, then

\be
K\oo p_\l(z_1,z_2)=c_l(q,t)\psi_\l(q,t;t^{3/2}z_1)\psi_\l(q,t;t^{3/2}z_2),
\label{A23}
\ee
where $z_1,\>z_2$ satisfy (\ref{A7}), the action of the operator
$K$ is given by formulas (\ref{A10}, \ref{A16}, \ref{A20}) and
the function $\psi_\l(q,t;z)$ is defined by (\ref{S8})}.

\vspace{0.3cm}
Using {\it MATHEMATICA} program we checked that factorization formula
(\ref{A23}) works for all partitions with $|\l|\le 4$.

Note that the separated equations (\ref{A6}) do not contain any information
about the normalization multiplier $c_\l(q,t)$ in (\ref{A23}).

\section{Discussion}

In this paper we constructed the discrete operator $K$ which should transform
$A_2$ Macdonald polynomials into the product of two $\phantom{|}_3\phi_2$
functions. It will be a theorem if conjecture (\ref{A22}) is
valid for the arbitrary integer $n$. In fact, formula (\ref{A22})
connects two very-well-poised $\phantom{|}_6\psi_6$ series with
different arguments. We failed to find relation (\ref{A22}) in any
literature related to basic hypergeometric series. So the proof of
(\ref{A22}) could be quite interesting from a purely mathematical point
of view.

The next problem is to calculate the inverse operator $K^{-1}$. We know
the answer in the limit $q\to1$ (see \cite{KS}) when $K$ is reduced
to the integral operator. But the problem of the calculation $K^{-1}$
in a discrete case can be much more difficult. Possibly the best way
to do this is to write a set of difference equations on the inverse
operator $K^{-1}$ and to solve them. Then we will obtain
a representation of $A_2$ Macdonald polynomials in terms of
basic $\phantom{|}_3\phi_2$ hypergeometric series.

But definitely a normalization of this discrete representation
will differ from the standard one (\ref{M6}). So a separate
problem is to calculate the corresponding
normalization factor.

And, of course, it would be interesting to generalize these results
for $A_{n-1}$ case. Separated equations (\ref{S10}) show that it can be
done. But here we come to another  problem:
what is an underlying algebraic structure for finite $q$ which
permits a separation of variables for Macdonald polynomials ?
For the classical case of Calogero-Sutherland
system the answer is known (see \cite{SK}).
So the answer can give a new insight on the whole algebraic
construction and help to write a discrete representation
for $A_n$ Macdonald polynomials. Possibly it  can give a lot of new
interesting identities in the theory of $q$-series.
We hope to consider  these problems in a separate paper.

\section{Acknowledgments}
I would like to thank Prof. G. Andrews
who explained me that the expression (\ref{A18}) at $n=0$
can be represented in the form of $\phantom{|}_6\psi_6$ series
and is summable by Bailey's formula.

\appendix

\section{Appendix}

In this appendix we write out  difference equations on the kernel of
the linear operator $K$ which appear from separated equations
(\ref{A6}).

Using explicit formulas for Macdonald operators (\ref{A1}) and
separated equations (\ref{A6}) we obtain four difference
equations on the kernel $K$:

\beq
&\Bigl\{(t-z_i)+(1-z_i{\ds q^2\over\ds t^3}){\ds(1-y_1)(1-y_2)\over
\ds
(1-ty_1)(1-ty_2)}T_{z_i}^2T_{y_1}T_{y_2}-(1-z_i{\ds
q\over\ds t^2})T_{z_i}\times&\nonumber\\
\label{Ap1}\\
&\times\Bigl[{\ds(ty_1-qy_2)(1-y_2)\over
\ds(y_1-qy_2)(1-ty_2)}T_{y_2}+{\ds(ty_2-qy_1)(1-y_1)\over
\ds(y_2-qy_1)(1-ty_1)}T_{y_1}\Bigl]
\Bigl\}K(z;y)=0,\nonumber&
\eeq
\beq
&\Bigl\{(1-z_i)t+(1-z_i{\ds q^2\over\ds t^2}){\ds(1-qy_1/t)(1-qy_2/t)\over
\ds(1-qy_1)(1-qy_2)}T_{z_i}^2T_{y_1}T_{y_2}-(1-z_i{\ds q\over\ds t})
T_{z_i}\times&
\nonumber\\
\label{Ap2}\\
&\times\Bigl[{\ds(ty_1-qy_2)(1-qy_2/t)\over
\ds(y_1-qy_2)(1-qy_2)}T_{y_2}+{\ds(ty_2-qy_1)(1-qy_1/t)\over
\ds(y_2-qy_1)(1-qy_1)}T_{y_1}\Bigl]\Bigl\}K(z;y)=0,&\nonumber
\eeq
where $i=1,2$.

Combining equations (\ref{Ap1}-\ref{Ap2}) one can obtain a
nonlinear difference equation in one variable, but  its solution
should satisfy some additional compatibility conditions.
In this way we obtained difference equations (\ref{UR1}-\ref{UR2})
which satisfy all necessary conditions.

After some tedious calculations one can check that the solution
(\ref{A16}) satisfies all equations (\ref{Ap1}-\ref{Ap2}) as well.

\section{Appendix}

In this appendix we calculate the normalization multiplier for the kernel $K$
and give three first nontrivial examples which support the conjecture
(\ref{A22}).

Using usual definitions for $\phantom{|}_r\psi_r$ series
(see \cite{GR,An}) we can rewrite
(\ref{A18}) as a sum of two $\phantom{|}_6\psi_6$ series:

\beq
\Phi(s,\sigma,\xi,t;q,n)={\ds q^s(1-q^{2s})\over\ds t^{3s}}
{\ds\Bigl({\ds q\sigma\over \ds t}q^{s},{\ds q\over\ds t\sigma}q^{s},
{\ds q\over\ds\xi\sqrt{t}}q^{s},{\ds q\xi\over\ds\sqrt{t}}q^{s}
\Bigl)_\infty
\over
{\ds\Bigl({\ds t\over \ds \sigma}q^{s},{t\sigma}q^{s},
{\xi\sqrt{t}}q^{s},{\ds \sqrt{t}\over\ds\xi}q^{s}
\Bigl)_\infty}}\times\phantom{aaa}&&\nonumber\\
\times\Biggl\{q^{ns}\phantom{|}_6\psi_6\Biggl[
\begin{array}{cccccc}
q^{s+1},& -q^{s+1},& {\ds t\over \ds\sigma}q^s,&
t\sigma q^s, & \xi\sqrt{t}q^s,& {\ds\sqrt{t}\over\ds\xi}q^s\\
q^{s},& -q^{s},& {\ds q\sigma\over \ds t}q^s,&
{\ds q\over \ds t\sigma} q^s,
& {\ds q\over \ds \xi\sqrt{t}}q^s,&
{\ds q\xi\over\ds\sqrt{t}}q^s\end{array};q,{\ds q^{1+n}\over\ds t^3}
\Biggl]+\phantom{a}&&\nonumber\\
+q^{-ns}\phantom{|}_6\psi_6\Biggl[
\begin{array}{cccccc}
q^{s+1},& -q^{s+1},& {\ds t\over \ds\sigma}q^s,&
t\sigma q^s, & \xi\sqrt{t}q^s,& {\ds\sqrt{t}\over\ds\xi}q^s\\
q^{s},& -q^{s},& {\ds q\sigma\over \ds t}q^s,&
{\ds q\over \ds t\sigma} q^s,
& {\ds q\over \ds \xi\sqrt{t}}q^s,&
{\ds q\xi\over\ds\sqrt{t}}q^s\end{array};q,{\ds q^{1-n}\over\ds t^3}
\Biggl]\Biggl\}.\>\>\label{B1}
&&
\eeq

Now we will calculate explicitly the function $\Phi(s,\sigma,\xi,t;q,0)$ using
Bailey's summation formula \cite{Ba} for the very-well-poised series
$\phantom{|}_6\psi_6$. One can check that  arguments
of $\phantom{|}_6\psi_6$ in (\ref{B1}) satisfy all necessary constraints
(at $n=0$) and so $\Phi(s,\sigma,\xi,t;q,0)$ is
summable. In fact, at $n=0$ we have a particular
case of Bailey's formula (which contains 5 independent parameters)
with 4 independent parameters.

So we have (see \cite{Ba,GR})
\beq
&\Phi(s,\sigma,\xi,t;q,0)=
{\ds2q^s(1-q^{2s})\bigl(q,q^{1-2s},q^{1+2s},q/t,q/t^2
\bigl)_\infty\over
\ds t^{3s}\bigl(q/t^3\bigl)_\infty}\times&\nonumber\\
&
\times{
\ds\Bigl({\ds q\sigma\xi\over \ds t^{3/2}},
 {\ds q\sigma\over\ds\xi t^{3/2}},{\ds q\xi\over\ds \sigma t^{3/2}},
{\ds q\over \ds\sigma\xi t^{3/2}}\Bigl)_\infty\over \ds\Bigl(
{\ds t\over\ds\sigma}q^s,t\sigma q^s,\xi\sqrt{t}q^s,
{\ds\sqrt{t}\over\ds\xi}q^s,
{\ds q\sigma\over \ds tq^s},{\ds q\over \ds t\sigma q^s},
{\ds q\over \ds\xi\sqrt{t}q^s},{\ds q\xi\over\ds \sqrt{t}q^s}\Bigl)_\infty
}.\label{B2}
&
\eeq

Now we give three first nontrivial examples of polynomials
$P(\sigma,\xi,t;q,n)$ (see (\ref{A22})):
\be
P(\sigma,\xi,t;q,1)={\ds t(\sigma+\sigma^{-1})+(t+1)\sqrt{t}(\xi+\xi^{-1})\over
\ds 2(1+t+t^2)},\label{B3}
\ee
\beq
&P(\sigma,\xi,t;q,2)={\ds
t(t(qt-1)(\sigma^2+\sigma^{-2})+(1+t)(qt^2-1)(\xi^2+\xi^{-2}))\over\ds
2(qt^3-1)(1+t+t^2)}+&\nonumber\\
&+{\ds(1+q)(t^2-1)
(t\sqrt{t}(\sigma+\sigma^{-1})(\xi+\xi^{-1})-(1+t)(1-t+t^2))\over\ds
2(qt^3-1)(1+t+t^2)},\label{B4}
&
\eeq

\beq
&P(\sigma,\xi,t;q,2)=
{\ds (1+q+q^2)(t^2-1)\sqrt{t}\over\ds 2(1+t+t^2)(q^2t^3-1)(qt^3-1)}
\times&\nonumber\\
&\Bigl[
\sqrt{t}(\sigma+\sigma^{-1})(1+(1+q)(t^3-t)-qt^4+t(qt^2-1)(\xi^2+\xi^{-2}))
+&\nonumber\\
&+(\xi+\xi^{-1})(1+(1+q)(t^3-t^2)-qt^5+t^2(qt-1)(\sigma^2+\sigma^{-2}))
\Bigl]+&\label{B5}\\
&
{\ds
t^3(qt-1)(q^2t-1)(\sigma^3+\sigma^{-3})+(1+t)(q^2t^2-1)(qt^2-1)t^{3/2}(\xi^3+\xi^{-3})
\over \ds 2(1+t+t^2)(q^2t^3-1)(qt^3-1)}.&\nonumber
\eeq

Apparently the function $P(\sigma,\xi,t;q,n)$ can be written in the form
of truncated basic hypergeometric series, but up to now
we failed to find this representation for the arbitrary $n$.

Also one can ask: is it possible to generalize formula (\ref{A22})
for the general 5-parametric case ? If the answer is yes,
then Bailey's summation formula
for $\phantom{|}_6\psi_6$ will be the first representative (the case $n=0$)
in a sequence of transformations for
very-well-poised $\phantom{|}_6\psi_6$ series
with different arguments.

\end{document}